\documentclass[12pt]{amsart}

\usepackage {amsmath}
\usepackage[applemac]{inputenc}
\usepackage {amsfonts, amssymb}
\usepackage {latexsym}

\begin{document}

\newtheorem {theorem} {Theorem} [section]
\newtheorem{lemma} [theorem]{Lemma} 
\newtheorem {cor}[theorem] {Corollary} 
\newtheorem {prop}[theorem] {Proposition} 
\newtheorem{preremark}[theorem]{Remark}
\newenvironment{rem}%
  {\begin{preremark}\rm}{\end{preremark}}
   \newtheorem{preremark1}[theorem]{Example}
\newenvironment{example}
  {\begin{preremark1}\rm}{\end{preremark1}}
  \newtheorem{preremark2}[theorem]{Definition}
\newenvironment{defn}%
  {\begin{preremark2}\rm}{\end{preremark2}}
 \renewcommand{\theequation}{\arabic{section}.\arabic{equation}}

\title[Asymptotic behaviour of 3D NS]{Asymptotic behaviour of the non-autonomous 3D Navier-Stokes problem with coercive force}

\author[D.Vorotnikov]{Dmitry Vorotnikov}

\address{CMUC, Department of
Mathematics, University of Coimbra, 3001-454 Coimbra, Portugal}{}
\email{mitvorot@mat.uc.pt}

\thanks{The research was partially supported by CMUC/FCT}




\begin{abstract}

We construct pullback attractors to the weak solutions of the three-dimensional Dirichlet problem for the incompressible Navier-Stokes equations in the case when the external force may become unbounded as time goes to plus or minus infinity. 

\end{abstract}

\maketitle

\section {Introduction} An \emph{attractor} of a dynamical system is a certain set to which every orbit eventually becomes close. When an autonomous differential equation (or boundary value problem) generates a dynamical system,  the corresponding attractor characterizes the long-time behaviour of its solutions \cite{hal,bv,cvb,har,tem}. The study of attractors to the 2D Navier-Stokes equations goes back to Ladyzhenskaya \cite{lad}, who was followed by lots of authors \cite{bab}. 

The non-autonomous equations do not automatically produce dynamical systems. Instead, one may define an attractor for a \emph{process} (a two-parameter semigroup) related to the solutions of a non-autonomous equation. There are three adequate approaches to this task. The first one is to extend the phase space and to deal with the \emph{skew-product} dynamical system \cite{se67,hal}. The second one \cite{cvb} is to introduce a concept of a \emph{uniform attractor} which attracts the trajectories uniformly with respect to the time shifts.  It turns out that sufficient conditions for existence of a uniform attractor \cite{cvb,c93} guarantee non-emptiness of the set (which is called the \emph{kernel} of the process) of bounded complete trajectories of the process. The \emph{sections} of the kernel possess \cite{c93} attraction properties which resemble the ones of the usual attractor of an autonomous system. However, this attraction is not uniform but \emph{pullback}, i.e. it happens when the actual moment of time is fixed and the initial time goes to minus infinity. The pullback mechanism appeared in much earlier works, e.g. \cite{kras} (see discussions in \cite{langa,gru}), but the concept of a \emph{pullback attractor} was proposed by Schmallfuss, Crauel and Flandoli (see \cite{or1,or2} and the references therein) in the early 1990s, and was since then successfully applied to many systems, e.g. \cite{min,cr,krs,langa,rrv,caidi1}.  This approach turned out to be relevant in much more general situations than the one of \cite{c93}\footnote{But the family of kernel sections, in the framework of \cite{c93}, coincides with the pullback attractor.}. Namely, the pullback attractors characterize the behaviour of processes rather at each ''finite'', ''present'' moment than as time goes to infinity. Therefore, this notion can be used to investigate the limiting behaviour of the processes which do not have bounded complete trajectories. Such a situation can arise, for instance, for equations with \emph{coercive}, i.e. unbounded as time goes to plus or minus infinity, right-hand members. The pullback attractors for the 2D Navier-Stokes system with (possibly) coercive non-autonomous body force were constructed in \cite{clr1,clr2}.

The attractor theory turned out to be generalizable onto the case of the problems which lack the property of uniqueness of solutions (or where such a property remains an open problem). Obviously, such problems do not generate dynamical systems in a normal manner. One of the main motivations for the progress in this direction was the ambition to study the limiting behaviour of the weak solutions to the 3D Navier-Stokes problem. There exist several ways of construction of attractors in this case. The first one, based on the theory of multivalued semigroups, goes back to \cite{bvi}, and was developed in \cite{mv}. It was used for the weak solutions of the 3D Navier-Stokes problem when the right-hand side is uniformly bounded in $H$ or under an unproved hypothesis \cite{kv}. A related \emph{generalised semiflow} approach was proposed in \cite{ball}, and adapted to stochastic problems in \cite{rubi}. 

An alternative method employs the concept of \emph{trajectory attractor}, i.e. the attractor of the translation semigroup in the space of trajectories \cite{sell,pures}. The sections of the trajectory attractor coincide with a properly defined \emph{global attractor} \cite{vc,cvb}. A similar procedure can be realized in the non-autonomous case at the presence of bounded complete trajectories, generalizing the notions of the uniform attractor and of the kernel \cite{pures,vc}.  The trajectory attractor technique is applicable to the weak solutions of the 3D Navier-Stokes problem \cite{sell,pures,vc,cvb}.  However, it requires the uniform boundedness of the Steklov average (in time) of the square of the $V^{*}$-norm of the body force. 

The trajectory attractor theory was amplified in \cite{book}, where some technical requirements, e.g. the invariance of the trajectory space, were omitted, which allowed us to study some problems where the classical trajectory attractor procedure was not working \cite{poly,vz1,vz2}. 

In \cite{cc}, the attractors to the 3D Navier-Stokes problem were handled in the framework of non-standard analysis. 

The treatments of pullback attractors for the non-autonomous problems without uniqueness are predominantly based on the concept of \emph{multi-valued dynamical
process} \cite{clmv,cmv,cmrv,wang,mareal,dcdss}. A trajectory attractor approach was introduced in
\cite{fs,fs1}, and, in a different manner, in \cite{caidi2}. The framework of \cite{caidi2} does not admit any unbounded trajectories.  The considerations of \cite{fs,fs1} were mainly directed at the analysis of stochastic equations; nevertheless, in \cite{fs1}, the deterministic 3D Navier-Stokes problem with unbounded body force was also studied. However, the coercivity was restricted by a complicated condition assuming some ''generalized boundedness'' as time goes to minus infinity (cf. \cite[p. 375]{fs1}), 
and the differentiability of the non-autonomous part of the forcing term in the spatial variable was supposed. 

In this work, we adapt the ideas from \cite{book} to the pullback attraction case. We introduce the notions of \emph{minimal pullback trajectory $\mathcal{D}$-attractor} and \emph{minimal pullback 
$\mathcal{D}$-attractor} (note that the latter is not a ''trajectory'' one). We find some general criteria for existence of these attractors. Then we investigate the relation between our concept of the minimal pullback 
$\mathcal{D}$-attractor and the existing one of the pullback $\mathcal{D}$-attractor for a dynamical process. Finally, we apply this approach to the construction of pullback attractors to the 3D Navier-Stokes problem. The only assumption on the body force is, roughly speaking, that the growth of its $V^{*}$-norm at minus infinity can be at most exponential. The same condition was imposed in \cite{clr1,clr2} for the 2D model.  

The paper is organized as follows. The next section is a preliminary one (notation etc.). The third section is devoted to the general description of our approach to the pullback attractors for the non-autonomous problems without uniqueness. The main results of the section are collected in Subsection 3.2, and the comparison with the pullback $\mathcal{D}$-attractor for a dynamical process
is carried out in Subsection 3.4. In the last section, we construct the minimal pullback trajectory $\mathcal{D}$-attractor and the minimal pullback 
$\mathcal{D}$-attractor for the weak solutions of the three-dimensional incompressible Navier-Stokes problem.

\newcommand {\R} {\mathbb{R}}
\newcommand {\E} {\mathbf{E}}
\def\be{\begin{equation}}
\def\ee{\end{equation}}
\def\fr#1#2{\frac{\partial #1}{\partial #2}}

\section{Preliminaries and notation}

Let $\Omega$ be a bounded domain (i.e. an open set, with any kind of boundary) in $\mathbb{R}^3$.

We shall use the standard notations $L_p (\Omega) $, $W_p^{\beta}
(\Omega) $, $H^{\beta} (\Omega) =$ $W_2^{\beta} (\Omega) $, $H^{\beta}_0 (\Omega) = \stackrel
{\circ}{W}{}_2^{\beta} (\Omega) $ $(\beta > 0)$ for the Lebesgue
and Sobolev spaces. 

Parentheses denote the following bilinear form:
$$ (u, v) = \int\limits_\Omega (u (x), v (x)) _ {F} dx, $$ where $F$ is $\R$, $\R^3$ or $\R^9$ (the space of 
$3\times 3$ - matrices).

The Euclid norm in $\R^3$ is denoted as $ | \cdot | $. The symbol $ \| \cdot \| $ will stand for the Euclid norm in $L_2(\Omega)$, $L_2(\Omega) ^3$, or $L_2(\Omega) ^9$. We shall also use the notation $ \| v \| _ 1 =
\| \nabla v \| $, $v\in H^1(\Omega)^3$.

Let $\mathcal {V}$ be the set
of smooth, divergence-free, compactly supported  in $ \Omega $  functions with
values in $\R^3$. The symbols $H $, $V $, $V_\delta $ ($ \delta>0 $)
denote the closures of $ \mathcal {V} $ in $L_2 (\Omega)^3$, $H^1(\Omega)^3, $ $H ^ {\delta} (\Omega)^3$, respectively.

Since $\Omega$ is bounded, there exists $\lambda_1>0$ so that 
\be\label{lam}\lambda_1\|u\|^2\leq \|u\|_1^2, \ u\in V.\ee

Following \cite {tem}, we identify the space $H $ and its
conjugate space $H ^ * $. Therefore we have the embedding
$$ V_\delta\subset H\equiv H ^*\subset V ^ *_\delta. $$
The value of a functional from $V ^ *_\delta $ on an element from 
$V_\delta $ is denoted by brackets $ \langle\cdot,
\cdot\rangle $.
We consider $V$ to be equipped with the norm $\|\cdot\|_1$ and
$V^*$ to be equipped with the corresponding norm of a conjugate
space.

The symbols $C (\mathcal{J}; E) $, $C_w (\mathcal{J}; E) $, $L_2
(\mathcal{J}; E) $ etc. denote the spaces of continu\-ous, weakly
continuous, quadratically integrable etc. functions on an interval
$\mathcal{J}\subset \mathbb {R} $ with
values in a Banach space $E $. We recall that a function $u:
\mathcal{J} \rightarrow E$ is \textit{weakly continuous} if for
any linear continuous functional $g$ on $E$ the function  $g(
u(\cdot)): \mathcal{J}\to \mathbb{R}$ is continuous.
Let us also remind that a pre-norm in the Frechet space $C ([0, +
\infty); E) $ may be defined by the formula
$$ \| v \| _{C ([0, +\infty);E)} =\sum\limits _ {i=1} ^ {+ \infty} 2 ^ {-i} \frac {\|v \| _ {C ([0, i]; E)}} {1 + \| v \| _ {C ([0, i]; E)}}. $$

Finally, let us introduce a very trivial notion, which will be useful to simplify the language. 

\begin{defn}A \emph{brochette} over a set $\mathcal{Y}$ is a family of sets $B_t\subset \mathcal{Y}$ depending on a scalar parameter $t\in\R$. \end{defn} 

To put it differently, a brochette is a multimap $B:\R\multimap \mathcal{Y}$.

\begin{defn} \label{def} For two brochettes $B$ and $B^*$ over $\mathcal{Y}$, we define the intersection $B\cap B^*$ as the family of $(B\cap B^*)_t=B_t \cap B^*_t$, $t\in \R$. We say that $B$ is contained in $B^*$ and write $B \subset B^*$ provided $B_t\subset B^*_t$ for all $t\in\R$. \end{defn} 

\section{Pullback trajectory and global attractors}

\subsection {Basic definitions}

Let $E $ and $E_0 $ be Banach spaces, $E\subset E_0 $.  Consider
an abstract non-autonomous differential equation\footnote{The symbol "$=$"\ may be understood in any appropriate sense (e.g.
in the sense of some topological space containing both $E$ and
$R(A)$). The derivative "$'$"\ may also be considered in any
generalized sense. The nonlinear operator $A$ is arbitrary (it may
even be multi-valued, but in this case the symbol "$=$" must be
replaced by "$\subset$").} \begin{equation}\label{main} u ' (t) =A (t,u
(t)),\end{equation}
$$ u:\R\to E,\ A: D(A)\to R(A), D(A)=\R\times E_A,\ E_A\subset E.$$

We study the limiting behaviour of the solutions to \eqref{main}
which continuously depend on time in the topology of  $E_0$.

We denote $\mathcal{T}=C ([0, + \infty); E_0) \cap L_{\infty,loc} (0, + \infty; E)$. Hereafter it is supposed that the space $E $ is reflexive. 
Then,
by a well-known Lions-Magenes lemma, see e.g. \cite[Lemma 2.2.6]{book}, $$\mathcal{T} \subset C_w ([0, + \infty); E).$$ Hence, the values
of functions from $\mathcal{T}$ belong to $E $ at every time.

We shall use the translation (shift) operators $T (h)$,
$$ T (h) (u) (t) =u (t+h), $$
where $h\geq 0 $ for $u\in \mathcal{T},$ 
and $h\in\mathbb{R}$ for $u\in C (\R;
E_0) \cup L_{\infty,loc} (\R; E) $.


For every $\tau\in \R$, let us consider some set
$$\mathcal {H} ^ +_\tau \subset \mathcal{T} $$
of solutions (strong, weak, etc.) to the shifted equation \begin{equation}\label{smain} u ' (t) =A (t+\tau,u
(t)),\end{equation} 
on the
positive time axis. The sets $ \mathcal {H} ^ +_\tau $ are
called \it trajectory spaces \rm and their elements are
called \it trajectories\rm. Note that $ \mathcal {H} ^ +$ is a brochette over $\mathcal{T}$ (the \emph{trajectory brochette}).


\begin{rem} An appropriate trajectory brochette  $ \mathcal {H} ^ +$ must be sufficiently ''wide'' in order to describe well the dynamics of \eqref{main}. Typically, it should be such that for every $a\in E$ and $\tau\in \R$ there exists (but is not necessarily
unique) a trajectory $u\in \mathcal {H} ^ +_\tau$ satisfying the initial condition
$u(0)=a$ (cf. \cite[Remark 4.2.2]{book} for the autonomous case).
\end{rem}
\begin{rem} As usual in the theory of trajectory attractors, the precise form of equation \eqref{main} is not significant (cf. \cite{cvb,book}). 
It merely matters to have a brochette $ \mathcal {H} ^ + $,
and everything depends on its properties only.
Generally speaking, the nature of $ \mathcal {H} ^ + $ may be
different from the one described above. \end{rem}

Now, fix a class $\mathcal{D}$ of brochettes $D=\{D_t\neq \varnothing,\ t\in \R\}$ over $E$. For each $D\in \mathcal{D}$, let us construct a brochette  $\mathcal{H}(D)$ according to the formula \be  \mathcal{H}_t (D)=\{u\in \mathcal{H}_t^+ : u(0)\in D_t\}.\ee

\begin{defn} A brochette $P$ over the set $\mathcal{T}$ is called \it pullback $\mathcal{D}$-attracting \rm (for $ \mathcal {H} ^ + $) if for all brochettes $ D\in\mathcal{D} $ and $t\in\R$ one has
$$\sup\limits _ {u\in  \mathcal{H}_\tau (D)} ^ {} \inf\limits _ {v\in P_t} ^ {} \|T (t-\tau) u-v \| _{C ([0, +\infty);E_0)} \underset {\tau\to-\infty} {\to} 0. $$

\begin{rem} This definition implies that, given a pullback $\mathcal{D}$-attracting brochette $P$, all the sets $P_t$ are non-empty.
\end{rem}

\end{defn} \begin{defn}
 A brochette $P$ over $\mathcal{T}$ is called  \it pullback $\mathcal{D}$-absorbing \rm (for 
$ \mathcal {H} ^ + $) if for all $ D\in\mathcal{D} $ and $t\in\R$ there is $\tau_0= \tau_0(D,t) \leq t $
such that for all $\tau \leq \tau_0 $ one has
$$ T (t-\tau)  \mathcal{H}_\tau (D)\subset P_t, $$ and the function $ \tau_0(D,\cdot):\R\to\R$ is non-decreasing for each fixed $D$. \end{defn}

It is easy to see that any absorbing brochette is an attracting one.

\begin{defn} A brochette $P$ over $\mathcal{T}$ is called  \it relatively $\mathcal{T}$-compact \rm if

i) $P_t$ is relatively compact in $C ([0, + \infty); E_0) $ for every $t\in\R$;

ii) there is a function $\phi:\R\times[0,\infty)\to\R,$ so that $\phi(t,\cdot)$ is continuous for fixed $t$, and $\|u(s)\|_E\leq \phi(t,s)$ for all $t\in \R, s\geq 0$ and $u\in P_t$.

Such a $P$ is called \it $\mathcal{T}$-compact \rm if, in addition,

i') $P_t$ is closed in $C ([0, + \infty); E_0) $ for every $t\in\R$.

\end{defn}

Given a brochette $P$ over $\mathcal{T}$, by $T (h) P$, $h\in\R$, we denote the following brochette: \be\label{th}(T (h) P)_t = T(h)P_{t-h}.\ee

\begin{defn} \label{SA} A brochette $P$ over $\mathcal{T}$ is called a \it pullback trajectory $\mathcal{D}$-semiattractor (PTSA) \rm for
$ \mathcal {H} ^ + $ if

i) $P $ is $\mathcal{T}$-compact;

ii) $T (h) P\subset P$ for any $h \geq 0 $ (in the sense of Definition \ref{def});

iii) $P $ is pullback $\mathcal{D}$-attracting.
\end{defn}

\begin{defn}  \label{PA} A PTSA is called a \it pullback trajectory $\mathcal{D}$-attractor (PTA) \rm for
$ \mathcal {H} ^ + $ if 

ii') $T (h) P = P$ for any $h \geq 0 $. \end{defn}

\begin{defn} A PTA is called a \it minimal pullback trajectory $\mathcal{D}$-attractor (MPTA) \rm for
$ \mathcal {H} ^ + $ if  it is contained (in the sense of Definition \ref{def}) in any other PTA.
 A PTSA is called a \it minimal pullback trajectory $\mathcal{D}$-semiattractor (MPTSA) \rm for
$ \mathcal {H} ^ + $ if  it is contained in any other PTSA.
\end{defn}

\begin{defn} \label{mpa} A brochette $\mathcal {A}$ over $E$ is called a \it minimal pullback 
$\mathcal{D}$-attractor (MPA) \rm for the trajectory brochette $ \mathcal {H}
^ + $  (in $E_0 $) if

i) $ \mathcal {A}_t $ is compact in $E_0 $ and bounded in $E $ for each $t\in\R$;

ii) for all $ D\in\mathcal{D} $ and $t\in\R$ there is pullback attraction:
$$\sup\limits _ {u\in  \mathcal{H}_\tau (D)} ^ {} \inf\limits _ {v\in \mathcal {A}_t} ^ {} \|u (t-\tau) -v \| _{E_0} \underset {\tau\to-\infty} {\to} 0. $$

iii) $ \mathcal {A} $ is the minimal brochette satisfying conditions i)
and ii) (i.e. $ \mathcal {A} $ is contained in every brochette
satisfying conditions i) and ii)). \end{defn}

\begin{rem} It is obvious that MPTA, MPTSA and MPA, if they exist, are unique. \end{rem}

\subsection {The main existence theorems}

\begin{theorem} \label{mt1} Assume that there exists a relatively $\mathcal{T}$-compact pullback $\mathcal{D}$-absorbing brochette $P $ for
$ \mathcal {H} ^ + $. Then there exists an
MPTA $\mathcal {U}\subset P$. 
\end{theorem}

\begin {theorem} \label{mt2} 
If a brochette $\mathcal{P}$ is a PTSA, then there exists an
MPTA $\mathcal {U}$ contained in $\mathcal{P}$. 
\end {theorem}

For a set $K\subset \mathcal{T}$, by $K(h)$, $h\geq 0$, we denote the set $\{v (h) |v\in K\}.$  Similarly, for a brochette $P$ over $\mathcal{T}$, by $P(h)$, $h\geq 0$, we denote the following brochette over $E$ (the \emph{section brochette}): $$(P(h))_t =  \{v (h) |v\in {P_t} \}.$$ 

\begin {theorem} \label{mt3}  If a brochette $ \mathcal {U} $ is an MPTA, then there is an MPA $ \mathcal
{A} $, and $\mathcal {A} = \mathcal {U} (0).$
\end {theorem}

\subsection {Proofs}

The proofs of the theorems require some preliminary observations.

\begin {lemma} \label{le4} a) For any two brochettes $P_1$ and $P_2$ over $\mathcal{T}$ satisfying the conditions
i) or ii) of Definition \ref{SA}, $P_1\cap P_2 $ also satisfies
a corresponding condition. b) If $P_1, P_2 $ are $\mathcal{T}$-compact and satisfy condition iii) of Definition
\ref{SA}, then $P_1\cap P_2 $ also satisfies condition iii).
\end {lemma}

\begin{proof} Statement a) is clear. Let us show b). Let $P_1$ and $P_2 $ be $\mathcal{T}$-compact and satisfy
condition iii). We have to show that $P_1\cap P_2 $ is a pullback 
$\mathcal{D}$-attracting set. If it is not so, then for some $ \delta> 0 $, $t\in\R$ and
$D\in \mathcal{D}$ there is a sequence $\tau_m\to -\infty $
such that
$$\sup\limits _ {u\in  \mathcal{H}_{\tau_m} (D)} ^ {} \inf\limits _ {v\in (P_1\cap P_2)_t} ^ {} \|T (t-\tau_m) u-v \| _{C ([0, +\infty);E_0)}> \delta. $$
Then there are elements $u_m\in  \mathcal{H}_{\tau_m} (D)$ such that \begin{equation} \label{con1}
\inf\limits _ {v\in (P_1 \cap P_2)_t} ^ {} \|T (t-\tau_m) u_m-v \| _{C ([0,
+\infty);E_0)}> \delta.\end{equation} On the other hand, since
$P_1 $ and $P_2 $ are pullback attracting, for any natural number $k $
there exist a number $m_k $ and elements $v^1_k\in (P_1)_t $,
$v^2_k\in (P_2)_t $ such that
$$ \| T (t-\tau_{m_k}) u _ {m_k}-v^1_k \| _{C ([0, +\infty);E_0)} <\frac 1 k, $$ $$ \| T (t-\tau_{m_k}) u _ {m_k}-v^2_k \| _{C ([0, +\infty);E_0)} <\frac 1 k. $$
Since $(P_1)_t $ is compact in $C ([0, + \infty); E_0) $, without loss
of generality we may assume that the sequence $v^1_k $ converges
to an element $v_0 $ as $k \to \infty $. Then the sequences $T (t-\tau_{m_k}) u _ {m_k} $ and $v^2_k $ also converge to $v_0 $. Thus, $
v_0\in (P_1\cap P_2)_t$ and $ \|T (t-\tau_{m_k}) u _ {m_k}-v_0 \| _{C
([0, +\infty);E_0)} \underset {k \to \infty} {\to} 0, $ which
contradicts \eqref{con1}. \end{proof} 

\begin {lemma} \label{le5} Let a brochette $P$ over $\mathcal{T}$ satisfy one of conditions
i), ii), ii') or iii) of Definitions \ref {SA} and \ref{PA}. Then $T (h) P $
also satisfies a corresponding condition for all $h \geq 0 $.
\end {lemma}

\begin{proof} Let $P $ satisfy condition ii), that is
$T(s)P_{t-s}\subset P_t$ for any $s \geq 0 $ and $t\in\R$. Then $$T(s)(T(h)P)_{t-s}= T(s)T(h)P_{t-s-h}$$ $$=T(h)T(s)P_{t-h-s}\subset T(h)P_{t-h}=(T(h)P)_t,$$ i.e. $T(h)P $ satisfies condition ii). The proof of the statement of the lemma concerning condition
ii') is similar, whereas concerning i) it is straightforward. Let $P $ satisfy condition iii), that is
it is pullback attracting. Since the map $T(h)$ is bounded in $C ([0, +
\infty ); E_0)$, one has $$\|T(h)u\|_{C ([0, + \infty ); E_0)}\leq
C \|u\|_{C ([0, + \infty ); E_0)}$$ for some constant $C$ and all
$u\in C ([0, + \infty ); E_0)$. Then for any $D\in \mathcal{D}$ and $t\in \R$ one has
$$\sup\limits _ {u\in  \mathcal{H}_{\tau} (D)} ^ {} \inf\limits _ {v\in T (h) P_{t-h}} ^ {} \|T (t-\tau) u-v \| _{C ([0, +\infty);E_0)}
= $$
$$\sup\limits _ {u\in  \mathcal{H}_{\tau} (D)} ^ {} \inf\limits _ {v\in P_{t-h}} ^ {} \|T (h) (T (t-h-\tau) u-v) \| _{C ([0, +\infty);E_0)}
\leq $$ $$ C\sup\limits _ {u\in  \mathcal{H}_{\tau} (D)} ^ {} \inf\limits _ {v\in P_{t-h}} ^
{} \|T (t-h-\tau) u-v \| _{C ([0, +\infty);E_0)} \underset {\tau \to-\infty}
{\to} 0, $$ and, due to \eqref{th}, $T (h) P $ is pullback attracting.  \end{proof} 

\begin {lemma} \label{le7} An MPTSA is always an MPTA.
\end {lemma}

\begin{proof} Let $ \mathcal {U} $ be an MPTSA. By Lemma \ref{le5}, $T
(h) \mathcal {U} $ is a PTSA for all $h \geq 0$, therefore $ \mathcal {U} \subset T (h) \mathcal {U} $. Thus, $
\mathcal {U} $ satisfies condition
ii') from Definition \ref{PA}, so it is a PTA, and obviously a minimal one.
\end{proof}

\begin{rem} \label{le8} The inverse statement is also true, but is based on Theorem \ref{mt2}, which we are still going to prove;  an MPTA is always an MPTSA.
Really, let $ \mathcal {U} $ be an MPTA and let $\mathcal{P}$ be a PTSA. By Theorem \ref{mt2}, $\mathcal
{U}\subset \mathcal{P}$. Thus, $ \mathcal {U}$ is contained in any
PTSA, so it is an MPTSA. \end{rem}

\begin{lemma} \label{le0} Assume that there exists a relatively $\mathcal{T}$-compact pullback $\mathcal{D}$-absorb\-ing brochette $P $ for
$ \mathcal {H} ^ + $. Then there is a
PTSA $\mathcal {P}\subset P$. 
\end{lemma}

\begin{proof} 
For every $ D\in\mathcal{D} $, $t\in\R$ and $\tau \leq \tau_0(D,t) $ one has
$ T (t-\tau)  \mathcal{H}_\tau (D)\subset P_t.$
Fix a number $t\in\R$, and take the closure in $C ([0, + \infty); E_0) $ of the
set
$$P^0_t=\bigcup\limits_{D\in \mathcal{D}}^{}\bigcup\limits_{\tau \leq \tau_0(D,t)}^{}T(t-\tau)\mathcal{H}_\tau (D),$$
and denote it by $\mathcal{P}_t$. The resulting brochette $\mathcal{P}$
is contained in $P$, therefore it is $\mathcal{T}$-compact. It is clear that it
is pullback absorbing. It remains to show that $T(h)P^0_{t-h}\subset P^0_t$ for $h\geq 0$. Then the continuity of the shift operator $T(h)$ in $C ([0, +\infty);E_0)$ 
would imply $T(h)\mathcal{P}_{t-h}\subset \mathcal{P}_t$, i.e. $T(h)\mathcal{P}\subset \mathcal{P}$.  Since the function $\tau_0(D,t) $ is non-decreasing in $t$, we have $$\bigcup\limits_{D\in \mathcal{D}}^{}\bigcup\limits_{\tau \leq \tau_0(D,t-h)}^{}T(t-\tau)\mathcal{H}_\tau (D) \displaystyle\subset \bigcup\limits_{D\in \mathcal{D}}^{}\bigcup\limits_{\tau \leq \tau_0(D,t)}^{}T(t-\tau)\mathcal{H}_\tau (D).$$ But the first union is $T(h)P^0_{t-h}$, and the second one is $P^0_t$.\end{proof}

\begin{lemma}\label{le6} (see \cite[Lemma 4.2.6]{book}) Let $ (X, \rho) $ be a metric space and
$ \{K_\alpha \} _ {\alpha\in\Xi} $ be a system of non-empty
compact sets in $X $. Assume that for any $ \alpha_1, \alpha_2 \in
\Xi $ there is $ \alpha_3\in\Xi $ such that $K _ {\alpha_1} \cap K
_ {\alpha_2} =K _ {\alpha_3} $.  Then $K_0 =\bigcap\limits _
{\alpha\in\Xi} ^ {} K_\alpha\neq \varnothing $, and for any $
\epsilon> 0 $ there is $ \alpha_\epsilon\in\Xi $ such that for any
$y\in K _ {\alpha_\epsilon} $
$$\inf\limits _ {x\in K_0} ^ {} \rho (x, y) <\epsilon. $$
\end {lemma}

Now we can begin to prove the theorems. 

\begin{proof} (Theorems \ref{mt1} and \ref{mt2}) We need to prove Theorem \ref{mt2}, and Theorem \ref{mt1} would then follow from Lemma \ref{le0}. 

Consider the intersection\footnote{Definition \ref{def} may evidently be generalized for the case of infinite number of intersecting brochettes.} $ \mathcal {U} $  of all pullback
trajectory $\mathcal{D}$-semiat\-tractors for
$ \mathcal {H} ^ + $.  
Let us show that $ \mathcal {U} $ is a PTSA.
Clearly, $ \mathcal {U} $ satisfies conditions i) and ii) of
Definition \ref{SA}. We are going to show that $ \mathcal {U} $
satisfies condition iii), i.e. it is pullback attracting. 

Fix $
\epsilon> 0 $, $t_0\in\R$ and a brochette $D\in\mathcal{D}$. In Lemma \ref{le6}, take $X=C ([0, + \infty ); E_0) $,
and let $ \{K_\alpha \} _ {\alpha\in\Xi} $ be the system of all sets $\mathcal{P}_{t_0}$ such that $\mathcal{P}$ is a
PTSA for $ \mathcal {H}
^ + $. By Lemma \ref{le4}, an intersection of two PTSAs is a PTSA, so the intersection of any
two sets from the system $ \{K_\alpha \} $ belongs to this system. It is clear that $$ \mathcal {U}_{t_0}=\bigcap\limits _
{\alpha\in\Xi} ^ {} K_\alpha.$$
By Lemma \ref{le6}, there is a
PTSA $\mathcal{P}_\epsilon $ such that for any $v\in (\mathcal{P} _ {\epsilon})_{t_0}
$
$$\inf\limits _ {w\in \mathcal {U}_{t_0}} ^ {} \|w-v \| _{C ([0, +\infty);E_0)} <\frac {\epsilon} 2. $$

Since $\mathcal{P}_\epsilon $ is a pullback attracting brochette, there exists $\tau_0$
such that, for $\tau\leq \tau_0$,
$$\sup\limits _ {u\in  \mathcal{H}_{\tau} (D)} ^ {} \inf\limits _ {v\in {(\mathcal{P}_\epsilon)_{t_0}}} ^ {} \|T (t_0-\tau) u-v \| _{C ([0, +\infty);E_0)} <\frac {\epsilon} 2. $$

Therefore for every $u\in  \mathcal{H}_{\tau} (D)$ there exists $v (u) \in (\mathcal{P}_\epsilon)_{t_0} $
so that
$$ \|T (t_0-\tau) u-v(u) \| _{C ([0, +\infty);E_0)} <\frac {\epsilon} 2. $$

We have: $$\sup\limits _ {u\in  \mathcal{H}_{\tau} (D)} ^ {} \inf\limits _ {w\in
\mathcal {U}_{t_0}} ^ {} \|T (t_0-\tau) u-w \| _{C ([0, +\infty);E_0)} \leq $$
$$\sup\limits _ {u\in  \mathcal{H}_{\tau} (D)} ^ {} (\|T (t_0-\tau) u-v (u) \| _{C ([0, +\infty);E_0)} + \inf\limits _ {w\in
\mathcal {U}_{t_0}} ^ {} \|v (u)-w \| _{C ([0, +\infty);E_0)}) $$ $$\leq
\frac {\epsilon} 2 +\frac {\epsilon} 2 =\epsilon. $$

Thus, $ \mathcal {U} $ is a PTSA, being the minimal one. By Lemma \ref{le7}, $ \mathcal {U} $ is an MPTA.
\end{proof}

\begin{proof} (Theorem \ref{mt3}) 
Observe first that the invariance property $T (h) \mathcal {U} =
\mathcal {U}, $ $ h\geq 0, $ implies $T (h) \mathcal {U}_{t-h} =
\mathcal {U}_t, $ and \be\label{for} \mathcal{U}_{t-h}(h) =\mathcal{A}_t \ee for every $t\in\R$, where $\mathcal{A}= \mathcal {U}(0)$. 

Every set $ \mathcal
{A}_t = \mathcal {U}_t (0) $, $t\in\R$, is  compact in $E_0 $ and bounded in $E $ due to $\mathcal{T}$-compactness of $\mathcal{U}$.

Take $ D\in\mathcal{D} $ and $t\in\R$. Since $ \mathcal {U} $ is a pullback attracting brochette,
$$\sup\limits _ {u\in  \mathcal{H}_\tau (D)} ^ {} \inf\limits _ {v\in \mathcal{U}_t} ^ {} \|T (t-\tau) u-v \| _{C ([0, +\infty);E_0)} \underset {\tau\to-\infty} {\to} 0. $$ It yields the pointwise
convergence:
$$\sup\limits _ {u\in  \mathcal{H}_\tau (D)} ^ {} \inf\limits _ {v\in \mathcal{U}_t} ^ {} \|(T (t-\tau) u-v)(h)\| _{E_0} \underset {\tau\to-\infty} {\to} 0,\ h\geq 0. $$ At $h=0 $ we get
$$\sup\limits _ {u\in  \mathcal{H}_\tau (D)} ^ {} \inf\limits _ {v\in \mathcal{A}_t} ^ {} \|u(t-\tau)-v\| _{E_0} \underset {\tau\to-\infty} {\to} 0. $$

It remains to show that $ \mathcal {A} $ is contained in every brochette
$A$ over $E$ satisfying the property 
\be \label{atta}\sup\limits _ {u\in  \mathcal{H}_\tau (D)} ^ {} \inf\limits _ {v\in A_t} ^ {} \|u(t-\tau)-v\| _{E_0} \underset {\tau\to-\infty} {\to} 0,   D\in\mathcal{D}, \ t\in\R,\ee and such that $A_t$ are compact in $E_0 $ and bounded in $E $.

Define a brochette $
U$ over $\mathcal{T}$ by the formula \be\label{ums}U_t= \{u\in \mathcal {U}_t |u (h) \in A_{t+h}
\forall h\geq 0 \}.\ee It suffices to show that $ \mathcal {U}
\subset U$. By Remark \ref{le8}, $ \mathcal {U} $ is
contained in every PTSA. Hence, it is enough
to show that $U$ is a PTSA. 

For any sequence $ \{u_m \} \subset U_t $ converging in $C ([0, +
\infty); E_0) $, its
limit $u_0 $ belongs to the (closed in $C ([0, + \infty); E_0) $)
set $ \mathcal {U}_t $. The convergence in $C ([0, + \infty); E_0) $
yields the pointwise convergence: $u_m (h) \to u_0 (h) $ in $ E_0,$ $h\geq 0 $. Since $ A_{t+h} $ is compact in $E_0 $,
$u_0 (h) \in A_{t+h} $, $h \geq 0 $. Thus, each $U_t$ is closed in $C ([0, + \infty); E_0) $. Since $
U\subset\mathcal {U} $, $U$ is
$\mathcal{T}$-compact.  Representation \eqref{ums} and the invariance property $T
(s) \mathcal {U} =\mathcal {U} $ yield $T (s) U \subset U,$ $s \geq 0 $. It remains to show
that $U$ is a pullback attracting brochette. 

If it is not so, then for some $ \delta> 0 $, $t\in\R$ and
$D\in \mathcal{D}$ there is a sequence $\tau_m\to -\infty $
such that
$$\sup\limits _ {u\in  \mathcal{H}_{\tau_m} (D)} ^ {} \inf\limits _ {v\in U_t} ^ {} \|T (t-\tau_m) u-v \| _{C ([0, +\infty);E_0)}> \delta. $$
Then there are elements $u_m\in  \mathcal{H}_{\tau_m} (D)$ such that \begin{equation} \label{co1}
\inf\limits _ {v\in U_t} ^ {} \|T (t-\tau_m) u_m-v \| _{C ([0,
+\infty);E_0)}> \delta.\end{equation} 

Since
$\mathcal{U}$ is pullback attracting, for any natural number $k $
there exist a number $m_k $ and elements $v_k\in \mathcal{U}_t $,
such that
$$ \| T (t-\tau_{m_k}) u _ {m_k}-v_k \| _{C ([0, +\infty);E_0)} <\frac 1 k.$$
 
But $ \mathcal{U}_t $ is compact in $C ([0, + \infty); E_0) $, so without loss
of generality we may assume that the sequence $v_k $ converges
to an element $v_0\in \mathcal{U}_t $ as $k \to \infty $. Then
\begin{equation} \label{co2} \|  T (t-\tau_{m_k}) u _ {m_k}-v_0 \| _{C ([0,
+\infty);E_0)} \underset {k \to \infty} {\to} 0. \end{equation}
Now \eqref{co1} and \eqref{co2} yield $v_0 \not\in U_t $, that is $v_0 (s) \not\in A_{t+s} $ for some $s
\geq 0 $. Using \eqref{atta} one gets
$$\inf\limits _ {v\in {A} _{t+s}} ^ {} \|u _ {m_k} (t+s-\tau_{m_k})-v \| _ {E_0} \underset {k\to\infty} {\to} 0. $$
Then there is a sequence $ \{v ^ * _ k \}\subset {A} _{t+s} $
such that
$$ \| T (t-\tau_{m_k}) u _ {m_k}(s)-v ^ * _ k \| _ {E_0} \underset {k\to\infty} {\to} 0. $$
Since $ {A} _{t+s} $ is compact, without loss of generality
$v_k ^ * $ converges to some element $v ^ * $. But  \eqref{co2} gives
$$ \|  T (t-\tau_{m_k}) u _ {m_k}(s)-v_0 (s) \| _ {E_0} \underset {k \to \infty} {\to}
0. $$ Therefore $v_0 (s) =v ^*\in {A} _{t+s} $, and we have
a contradiction. \end{proof}

\subsection {A comparison of the concept of MPA with the pullback $\mathcal{D}$-attractors for a process}


We keep assuming that we are given some spaces $E$, $E_0$ and a fixed class $\mathcal{D}$ of brochettes over $E$.
We recall that a \emph{process} $U$ on $E$ is a two-parameter family of maps $$U(t,\tau):E\to E,\ t,\tau\in \R,\ t\geq \tau,$$ so that 
$U(t,t)\xi=\xi$ and $U(t,\tau)\xi=U(t,s)U(s,\tau)\xi,$ for all $\xi \in E$ and $t,s,\tau\in\R, \ t\geq s\geq\tau$.

\begin{defn} \label{pro} A brochette $\mathcal {A}$ over $E$ is called a \it pullback 
$(E, E_0,\mathcal{D})$-attractor \rm for $U$ if

i) $ \mathcal {A}_t $ is compact in $E_0 $ and bounded in $E $ for each $t\in\R$;

ii) $\mathcal {A}$ is \it pullback 
$(E, E_0,\mathcal{D})$-attracting \rm for $U$, that is
\be\sup\limits _ {u\in D_\tau} ^ {} \inf\limits _ {v\in \mathcal {A}_t} ^ {} \|U(t,\tau)u-v \| _{E_0} \underset {\tau\to-\infty} {\to} 0\ee for all $D\in\mathcal{D}$ and $t\in\R$; 

iii) $\mathcal {A}$ is \emph{invariant}, i.e. \be U(t,\tau) \mathcal {A}_\tau = \mathcal {A}_t \ee for $t,\tau\in\R, \ t\geq\tau$. \end{defn}

\begin{rem} This definition is equivalent to a standard one (see e.g. \cite{clr1,clr2}) in the case $E=E_0$. For the sake of generality, we consider the general case $E\subset E_0$, where the topology of attraction (in our case, the one of $E_0$) may be different from the one of the phase space $E$ (see e.g. \cite{bv,cvb} for similar approaches to attractivity). \end{rem} 

\begin{rem} Pullback 
$(E, E_0,\mathcal{D})$-attractors, as defined above, can be not unique (a simple example may be found in \cite{min}). Some minimality conditions (see e.g. \cite{rrv,min}) may be added to the definition in order to provide uniqueness (we return to this issue below, in Remark \ref{minn}). \end{rem} 

Processes are usually generated by non-autonomous differential equations. Assume that for any $b\in E$ and $\tau\in\R$, equation \eqref{main} possesses a unique
solution $$u_{b,\tau}\in C ([\tau, +
\infty); E_0) \cap L_{\infty,loc} (\tau, + \infty; E),$$ satisfying the
initial condition \begin{equation} u_{b,\tau}(\tau)=b.\end{equation}

Then one can define the process $U$ corresponding to \eqref{main} by the formula
\begin{equation} \label{c1} U(t,\tau)(\xi)= u_{\xi,\tau}(t).\end{equation}

In this situation the natural family of
trajectory spaces is
\begin{equation} \label{c2} \mathcal {H}^+_\tau=\{u_b\in \mathcal{T}| u_b(\cdot)= U(\cdot+\tau,\tau) b, b\in
E\}, \tau\in\R.\end{equation} 

Now we examine the relation between Definitions \ref{pro} and \ref{mpa}.

\begin{theorem} \label{equ} a)  If there exists a pullback 
$(E, E_0,\mathcal{D})$-attractor $\mathcal{A}$ for $U$, and $\mathcal{A}\in \mathcal{D}$, then $\mathcal{A}$ is an
MPA for $ \mathcal {H} ^ +$. b) Let the conditions of Theorem
\ref{mt3} hold for the trajectory brochette $\mathcal {H}^ +$.  If the MPTA
$\mathcal{U}$ is contained in $\mathcal {H}^ +$ (in the sense of Definition \ref{def}), then the MPA $\mathcal{A}=\mathcal{U}(0)$ is a
pullback 
$(E, E_0,\mathcal{D})$-attractor for $U$.
\end{theorem}
\begin{proof}  Due to the identity
\begin{equation}\sup\limits _ {u\in  \mathcal{H}_\tau (D)} ^ {} \inf\limits _ {v\in \mathcal {A}_t} ^ {} \|u (t-\tau) -v \| _{E_0}  =\sup\limits _ {b\in D_\tau} ^ {} \inf\limits _ {v\in \mathcal {A}_t} ^ {} \|U(t,\tau)b-v \| _{E_0} 
\end{equation} for all $t,\tau\in \R$, $t\geq\tau$, and $D\in\mathcal{D}$, conditions i) (which simply coincide) and ii), resp.,  of Definitions \ref{mpa} and \ref{pro},  are equivalent\footnote{Of course, under assumptions \eqref{c1} and \eqref{c2}.}. To prove a),  it remains to show that a pullback 
$(E, E_0,\mathcal{D})$-attractor $ \mathcal {A}\in\mathcal{D}$ for $U$ is contained in any brochette $A$ for which axioms i) and ii) of Definition \ref{pro} hold. 
Fix an arbitrary number $t\in\R$. Since $A_t$ is compact in $E_0$, for any open neighborhood $W$ of $A_t$ in $E_0$ one has $U(t,\tau) \mathcal{A}_\tau \subset W$ for all
$\tau$ close to $-\infty$. If there is a point $w\in \mathcal{A}_t$ such
that $w\not\in A_t$, then $W_w=E_0\backslash\{w\}$ is an open
neighborhood of $A_t$. Therefore $w\in\mathcal{A}_t=U(t,\tau) \mathcal{A}_\tau\subset W_w$, and we
arrive at a contradiction. 

To check b), we only need to show that, under the conditions of Theorem
\ref{mt3}, the brochette  $\mathcal{A}=\mathcal{U}(0)$ is invariant. But the inclusion $\mathcal{U}\subset\mathcal {H} ^ +$ and representation
\eqref{c2} yield
\begin{equation}\mathcal{U}_\tau=\{u_b\in \mathcal{T}| u_b(\cdot)= U(\cdot+\tau,\tau) b, b\in
\mathcal{U}_\tau(0)\}, \tau\in\R.\end{equation} 

Hence, for all $t\geq \tau$, $$U(t,\tau)
\mathcal{U}_\tau(0)=\mathcal{U}_\tau(t-\tau),$$ and by \eqref{for} we conclude: $$U(t,\tau)
\mathcal{A}_\tau=\mathcal{A}_t.$$
 \end{proof}
 
\begin{rem} \label{minn} The above argument shows that a pullback 
$(E, E_0,\mathcal{D})$-attractor is in a certain sense minimal provided it belongs to the set $\mathcal{D}$. Note that the proof of this issue did not use the particular structure of the process $U$ and is thus valid for any process. Hence, the requirement for a pullback $(E, E_0,\mathcal{D})$-attractor to belong to $\mathcal{D}$ may be a relevant alternative to minimality constraints\footnote{By the way, an artificial a posteriori procedure can insure this condition. It suffices to replace $\mathcal{D}$ with $\mathcal{D}'=\mathcal{D}\cup \{\mathcal{A}\}$, where $\mathcal{A}$ is the given $(E, E_0,\mathcal{D})$-attractor. Then $\mathcal{A}$ is a $(E, E_0,\mathcal{D}')$-attractor belonging to the set $\mathcal{D}'$.}. For instance, the pullback attractors considered in \cite{clr1} meet this requirement. \end{rem} 

\section{Pullback attractors for the 3D Navier-Stokes problem}

\subsection{Weak solutions to the 3D Navier-Stokes problem}

Consider the 3D incompressible Navier-Stokes problem:

\begin{equation}\label{ns1} \frac {\partial u} {\partial t} + \sum\limits _ {i=1} ^ {3} u_i \frac {\partial
u} {\partial x_i}-\eta\Delta u +\nabla  p =F, \end{equation}

\begin{equation}\label{ns2} div\ u =
0,\end{equation}
\begin{equation}\label{ns3} u \Big | _ {\partial\Omega} =0,\end{equation}
where $u $ is an unknown velocity vector, $p $ is an unknown
pressure function, $F$ is the given body force (all of them
depend on a point $x $ in a bounded domain $ \Omega \subset \mathbb{R}^3$,
and on a moment of time $t$), and $ \eta>0 $ is the
viscosity of a fluid.
\begin{defn} Let $F\in L_{2,loc} (0, \infty; V^*)$. A function \be\label{cla} u\in L_{2,loc} (0, \infty; V)
\bigcap C_w ([0, \infty); H) \bigcap W^1_{4/3,loc} (0, \infty; V ^ *)
\ee
is an \textit{admissible weak solution} to problem \eqref{ns1}--\eqref{ns3}
if it is a weak solution, i.e.
\begin{equation}\label{weak}\frac {d} {d t} (u, \varphi) + \eta (\nabla u, \nabla \varphi)-\sum\limits^3 _ {i=1}
(u_i u, \frac {\partial\varphi} {\partial x_i}) = \langle F,
\varphi\rangle\end{equation}
for all test functions $ \varphi\in V $ a.e. on $ (0, \infty) $ (cf. e.g. \cite{temam}), and it satisfies the energy inequality \be\label{ener} \|u(h)\|^2\leq e^{-\sigma h}\left(\|u(0)\|^2+\frac 1 \eta \int\limits_{0}^h e^{\sigma \xi}\|F(\xi)\|^2_{ V^*}\,d\xi\right)\ee for all $h\geq 0$, where \be\label{seg}\sigma=\eta\lambda_1.\ee
\end{defn}

\begin{prop} For every $a\in H$ and $F\in L_{2,loc} (0, \infty; V^*)$, there exists an admissible weak solution to \eqref{ns1}--\eqref{ns3} satisfying the initial condition \begin{equation} u | _ {t=0} =a.\end{equation}\end{prop}

\begin{proof} Consider a family of approximating problems: find $$ u_M\in L_{2} (0, M; V)
\bigcap C ([0, M]; H) \bigcap W^1_{2} (0, M; V ^ *), \ u_M(0) =a,
$$ so that \begin{equation}\label{weak1} \langle u^\prime_M, \varphi \rangle + \eta (\nabla u_M, \nabla \varphi)-\sum\limits^3 _ {i=1}
(\frac{(u_M)_i u_M}{1+|u_M|^2/M}, \frac {\partial\varphi} {\partial x_i}) = \langle F,
\varphi\rangle\end{equation}
for all test functions $ \varphi\in V $ a.e. on $ (0, M) $, where $M$ is a natural number. It is known \cite{zvdm} that such problems possess solutions. 

We recall the identity (cf. \cite[p. 29]{aaa} or \cite[Formula (6.1.21)]{book}) \be\label{xx} \sum\limits^3 _ {i=1}
(\frac{u_i u}{1+|u|^2/M}, \frac {\partial u} {\partial x_i})=0, \ u\in V. \ee  
Substitute $ 2 e^{\sigma t}u_M(t)$ for $\varphi$ into \eqref{weak1} at a.a. $t\in(0,M)$: 

\be 2 e^{\sigma t} \langle u^\prime_M(t),  u_M(t) \rangle= - 2\eta e^{\sigma t}\| u_M(t)\|^2_1 +2 e^{\sigma t}  \langle F(t),
u_M(t)\rangle.\end{equation}

This implies

\be \frac {d} {d t} (e^{\sigma t} \|u_M(t)\|^2)  -\sigma e^{\sigma t}\| u_M(t)\|^2 $$ $$ \leq - \eta e^{\sigma t}\| u_M(t)\|^2_1  + \frac 1 \eta e^{\sigma t} \|F(t)\|^2_{ V^*}. \ee

Integrating from $0$ to $s\geq 0$, and taking into account \eqref{lam} and \eqref{seg}, we get
\be e^{\sigma s}\|u_M(s)\|^2\leq \|a\|^2+\frac 1 \eta \int\limits_{0}^s e^{\sigma \xi}\|F(\xi)\|^2_{ V^*}\,d\xi.\ee

Therefore, for all $h\geq 0$,
\be \label{inqc}\max_{0\leq s\leq h} e^{\sigma s}\|u_M(s)\|^2\leq \|a\|^2+\frac 1 \eta \int\limits_{0}^h e^{\sigma \xi}\|F(\xi)\|^2_{ V^*}\,d\xi.\ee

Due to \eqref{xx}, the solutions to \eqref{weak1} satisfy the standard bounds on $\|u_M(t)\|$ and $\int\limits_{0}^t \|u_M(\xi)\|^2_{1}\,d\xi$ available for the weak solutions of the Navier-Stokes problem \cite{lions,temam}, uniformly with respect to $M$.  Via a diagonal argument one easily concludes that there exist a subsequence $u_{M_k}$ and a limiting function $u$ such that $u_{M_k}\to u$ as $k\to\infty, M_k>T$, weakly in $L_2(0,T;V)$, weakly-* in $L_\infty(0,T;H)$, and strongly in $L_2(0,T;H)$ for every $T>0$. This function $u$ is a weak solution to \eqref{ns1}--\eqref{ns3} in class \eqref{cla}.
Passing to the limit in \eqref{inqc}, we get \be \label{inqc1}\mathrm{ess} \sup_{0\leq s\leq h} e^{\sigma s}\|u(s)\|^2\leq \|a\|^2+\frac 1 \eta \int\limits_{0}^h e^{\sigma \xi}\|F(\xi)\|^2_{ V^*}\,d\xi.\ee This yields
\be e^{\sigma h}\|u(h)\|^2\leq \|a\|^2+\frac 1 \eta \int\limits_{0}^h e^{\sigma \xi}\|F(\xi)\|^2_{ V^*}\,d\xi,\ee
e.g. by \cite[Theorem 1.7, p. 33]{cvb}.
\end{proof}

\subsection{Minimal pullback attractors for 3D NS}

Fix $f\in L_{2,loc} (\R; V^*)$ such that \be\int\limits_{-\infty}^t e^{\sigma \xi}\|f(\xi)\|^2_{ V^*}\,d\xi < +\infty\ee
for some (and thus for all) $t\in\R$.
Let us construct an MPTA and an MPA for the Navier-Stokes problem \eqref{ns1}--\eqref{ns3} with $F=f$. 

We take $$E=H,$$ and $$E_0
=V_\delta ^*,$$ where $\delta \in (0,1] $ is a fixed
number. We define $\mathcal{D}$ as follows (cf. \cite{clr1,clr2}). Let $\mathcal{R}$ be the set of such functions $r:\R\to (0,+\infty)$ that  \be\lim\limits_{s\to -\infty}^{} e^{\sigma s} r^2 (s)=0,\ee and the function $e^{\sigma \cdot} r^2$ is increasing. The class $\mathcal{D}$ consists of the brochettes $D$ over $H$ for which there exist functions $r_D\in \mathcal{R}$ so that $\|w\|\leq r_D(t)$ for all $t\in \R$ and $w\in D_t$. 


The trajectory spaces $ \mathcal {H}
^ +_\tau $, $\tau\in\R$, are the sets of admissible  weak solutions to \eqref{ns1}--\eqref{ns3} with the shifted right-hand members $F=T(\tau)f$. These trajectory spaces are contained in $\mathcal{T}$. In fact, by \eqref{ener},
every admissible weak solution $u$ 
belongs to $L_{\infty,loc} (0, + \infty; H)$. Since $ \Omega $ is bounded, $V_\delta\subset H $
compactly, thus $H\subset V ^
*_\delta $ compactly. But $u'\in L_{4/3,loc} (0, \infty; V ^ *)$, so $u\in C ([0,
\infty); V ^ *_\delta) $ by the Aubin-Simon compactness theorem \cite[Corollary 4]{sim}.

\begin {theorem} For the trajectory brochette
$ \mathcal {H} ^ + $, there exist an MPTA $ \mathcal {U}$ and an MPA $\mathcal {A} = \mathcal {U} (0)$. Moreover, $\mathcal {A} \in \mathcal {D}$.
\end {theorem}

\begin{proof}  

Consider the brochette $P$ over $\mathcal{T}$ so that the sets $P_t$, $t\in\R$, consist of functions $u\in \mathcal{T}$ satisfying
the inequalities
\be\|u(h) \|^2 \leq 2 e^{-\sigma (t+h)} R_1(t+h), \ee
\be \|u ' (h)\| _ {V_3^*} \leq \eta  R_2 \|u(h) \| +R_3 \|u(h) \|^2+\|f(t+h)\|_ {V_3^*} \ee
for a.a. $h\geq 0$, where $$R_1(s)=e^{\sigma s}+\frac {1} \eta \int\limits_{-\infty}^s e^{\sigma \xi}\|f(\xi)\|^2_{ V^*}\,d\xi,$$ and the constants $R_2$ and $R_3$, depending only on the domain $\Omega$, will be defined below.

By \cite[Corollary 4]{sim}, the sets $P_{t,M} =
\{v=u | _ {[0, M]}: u\in P_t \} $, $M>0$, are relatively compact in $C ([0,
M]; E_0) $. This immediately implies (cf. e.g. \cite[p. 183]{book})
that $P_t$ are relatively compact in $C ([0, + \infty); E_0) $. Now it is easy to conclude
that $P$ is relatively $\mathcal{T}$-compact.  

Let us check that the brochette $P $ is pullback $\mathcal {D}$-absorbing. Fix $t\in\R$ and $D\in\mathcal{D}$. 
Set $$\chi(s)=\max\{e^{\sigma s} r^2_D(s),\, R_1(s)\},\ s\in\R.$$ 
Note that the functions $R_1$ and $\chi$ are increasing. Thus, $\tau_0=\chi^{-1}(R_1(t))$ is an increasing 
 function of $t$ (for fixed $D$), and $\tau_0\leq t$. Let $\tau\leq \tau_0$. We have to show that $T (t-\tau)  \mathcal{H}_\tau (D)\subset P_t$. Let $u\in  \mathcal{H}_\tau (D)$, i.e. $u\in \mathcal{H}_\tau^+ $ and $u(0)\in D_\tau$. Due to \eqref{ener}, for the function $v=T (t-\tau)u$ and a.a. $h\geq 0$, we have
\be \|v(h)\|^2=\|u(t-\tau+h)\|^2 $$ $$\leq e^{-\sigma (t-\tau+h)}\left(\|u(0)\|^2+\frac 1 \eta \int\limits_{0}^{t-\tau+h}
 e^{\sigma \xi}\|f(\xi+\tau)\|^2_{ V^*}\, d\xi\right) $$ $$\leq e^{-\sigma (t+h-\tau)}r_D^2(\tau)+\frac 1 \eta
  \int\limits_{\tau}^{t+h} e^{\sigma (\xi-t-h)}\|f(\xi)\|^2_{V^*}\, d\xi $$ $$\leq e^{-\sigma (t+h)}\left[\chi(\tau)+R_1(t+h) \right]\leq 2e^{-\sigma (t+h)}R_1(t+h),\ee 
since $$\chi(\tau)\leq\chi(\tau_0)=R_1(t)\leq R_1(t+h).$$

The function $u$ satisfies \eqref{weak} with $F=T(\tau)f$, so $v$ satisfies \eqref{weak} with $F=T(t)f$. Take any function $\varphi\in V_3$. Then \be |\langle v'(h),\varphi \rangle|\leq \eta |(v(h), \Delta \varphi)|+\sum\limits^3 _ {i=1}
|(v_i(h) v(h), \frac {\partial\varphi} {\partial x_i})| + |\langle f(t+h),
\varphi\rangle| $$ $$\leq (\eta  R_2 \|v(h) \| +R_3 \|v(h) \|^2+\|f(t+h)\|_ {V_3^*})\|\varphi\|_{V_3},\ee with certain constants 
$R_2$ and $R_3$, depending only on the domain $\Omega$. We have applied the fact of the continuous Sobolev imbedding $$V_3\subset H^3_0(\Omega)^3\subset W^1_\infty(\Omega)^3$$ in 3D. 

Let $$r_A(t)=\sqrt{2 e^{-\sigma t} R_1(t)}.$$ Then $r_A\in \mathcal{R}$, and $\|w\|\leq r_A(t)$ for all $t\in \R$ and $w\in P_t(0)$. 

By Theorem \ref{mt1} there exists an MPTA $ \mathcal {U}$, and by Theorem \ref{mt3} there is an MPA $\mathcal {A} = \mathcal {U} (0)$. Finally, since $\mathcal {A}_t=\mathcal {U}_t (0)\subset P_t (0)$, we have $\mathcal {A} \in \mathcal{D}$.

\end{proof} 

\bibliography{brochette2}

\bibliographystyle{abbrv}
\end{document}